\documentclass[12pt]{amsart}
\usepackage{amsthm,amssymb}
 
\newtheorem{thm}{Theorem}[section]
\newtheorem{lm}[thm]{Lemma}

\newtheorem{pro}[thm]{Proposition}

\input epsf

\def \B {\mathcal B}

\def \o {\omega}

\begin{document}
 
\baselineskip.525cm

\title[Random walk, Jones polynomial and zeta function]
{Random Walk on knot diagrams, Colored Jones Polynomial and 
Ihara-Selberg Zeta Function}
\thanks{The first author is supported in part by an N.S.F. grant 
and the second author is supported by an NSF postdoctoral fellowship.}
\author[X.-S. Lin]{Xiao-Song Lin}
\address{Department of Mathematics, University of California, 
Riverside, CA 92521}
\email{xl@math.ucr.edu}
\author[Z. Wang]{Zhenghan Wang}
\address {Department of Mathematics, Indiana University, Bloomington, 
IN 47405}
\email{zhewang@indiana.edu}
\begin{abstract} 
A model of random walk on knot diagrams is used to study the 
Alexander polynomial and the colored Jones polynomial of knots.  
In this context, the inverse 
of the Alexander polynomial of a knot plays the role of an Ihara-Selberg 
zeta function of a directed weighted graph, counting with weights cycles of 
random walk on a 1-string link whose closure is the knot in question. 
The colored Jones polynomial then counts with weights families of 
``self-avoiding'' cycles of random walk on the cabling of the 1-string link. 
As a consequence of such interpretations of the Alexander and colored Jones
polynomials, the computation of the limit of the renormalized colored Jones polynomial 
when the coloring (or cabling) parameter tends to infinity whereas the weight parameter 
tends to 1 leads immediately to a new proof of the Melvin-Morton conjecture, which was 
first established by Rozansky and by Bar-Natan and Garoufalidis.
\end{abstract}
\maketitle

\section{Introduction}

\noindent The Alexander polynomial and the Jones polynomial, both 
characterized by simple crossing change formulae, are probably the two most 
celebrated invariants in knot theory.                                 
While the Alexander polynomial appears again and again in different contexts, 
making us feel quite comfortable with it, 
the nature of the Jones polynomial remains mysterious. 
In this paper, we will provide a new perspective for the study of the Jones
polynomial (and its generalizations -- the so-called colored Jones
polynomial), the Alexander polynomial and their relationship. An immediate
outcome of this new perspective is a straightforward proof of the
Melvin-Morton conjecture.

In \cite{LTW}, a model of random walk on knot diagrams was introduced. When 
we were seeking formulations of the Alexander and Jones polynomials in this 
model of random walk, a paper of Foata and Zeilberger \cite{FZ} caught 
our attention. In that paper, Foata and Zeilberger established a general
combinatorial framework for counting with weights Lyndon words in 
a free monoid generated by a totally ordered set, one of its consequences 
is a proof of Bass' evaluations of the Ihara-Selberg zeta function for 
graphs. We noticed that one of the main theorems of \cite{FZ} implies the
following fact: Take a 1-string link and 
consider all families of cycles on
this 1-string link in our model of random walk. Every cycle is assigned
with a weight (probability). Then the Ihara-Selberg type zeta function
constructed using these weights is equal to the inverse of the
Alexander polynomial of the knot obtained as the closure of the
1-string link, up to a factor in the form of a power of the weight parameter.

There is a remarkable relation between the colored Jones polynomial
and the Alexander polynomial, which was first noticed and conjectured by
Melvin and Morton \cite{MM}. Rozansky \cite{RO} gave an argument for 
this conjecture using the Chern-Simons path integral formalism of the
colored Jones polynomial and the relation between Ray-Singer analytic
torsion and the Alexander polynomial. The rigorous proof of the
Melvin-Morton conjecture was given by Bar-Natan and Garoufalidis \cite{BG},
using the full power of the theory of finite type knot invariants.

In our setting of random walk on knot diagrams, the Jones polynomial
counts only simple families of cycles on the 1-string link, 
i.e. families of cycles which do not share any edge. To take into account of
all cycles, we have to use the colored Jones polynomial. A state sum
formula for the (renormalized) colored Jones polynomial with the coloring 
parameter $d+1$
implies that it counts simple
families of cycles on $d$-cabling of the 1-string link in question. 
To relate the colored Jones polynomial with the Alexander polynomial, we lift 
families of cycles on the string link to its $d$-cabling with the weight 
parameter adjusted appropriately. A family of cycles on the 1-string link can 
have many liftings to its cabling. Weights of all liftings
add up to the weight of the original family of cycles,
whereas the weights of non-simple liftings vanish in the limit when 
$d\rightarrow\infty$. So in the limit, only weights 
of simple families of cycles survive and this calculation leads to a
proof of the Melvin-Morton conjecture.

We remark that our formulation of the limit of the colored Jones
polynomial is analogous to the limit of partition functions on a finite
lattice with a fixed boundary condition in statistical mechanics. Our 
proof of the Melvin-Morton conjecture is in spirit close to Rozansky's proof 
using the semi-classical limit of Chern-Simons path integral. 

The model of random walk on knot diagrams has a much richer content than we
have touched upon here. A more detailed exploration of this model
will be the subject of our future publications.

\section{Random walk on knot diagrams}
\subsection{Wirtinger presentation and free derivatives}

Fix an oriented knot diagram $K$, we will label the arcs in the knot diagram 
separated by crossings at the under-crossed strands using 
the letters $x_1,x_2,\dots, x_n$.  The knot group 
$G(K)=\pi_1({\mathbb R}^3\setminus K)$ admits a 
Wirtinger presentation 
as follows: It has $x_1,x_2,
\dots,x_n$ as generators, and one relation for each crossing. If a crossing 
has incident 
arcs $x_i,x_j,x_k$, where $x_i$ separates $x_j$ and $x_k$ in a small 
neighborhood of the crossing 
and the knot orientation points $x_j$ toward $x_k$, 
the relation is 
$$x_j=x_i^{\epsilon}x_kx_i^{-\epsilon}.$$
Here $\epsilon=\pm1$ is the sign of the crossing. 

With respect to the abelianization $\phi:{\mathbb Z}G(K)\rightarrow{\mathbb Z}
[t^{\pm1}]$, sending each 
$x_i$ to $t$, a free derivative $\partial:{\mathbb Z}G(K)\rightarrow{\mathbb
Z}[t^{\pm1}]$ is a linear map such that 
$$\partial(g_1g_2)=\partial(g_1)+\phi(g_1)\partial(g_2)\qquad\text{for
all $g_1,g_2\in G(K)$}.$$
The ${\mathbb Z}[t^{\pm1}]$-module of free derivatives on the free group 
$F$ generated by 
$x_1,x_2,\dots,x_n$ is spanned by 
$\partial_i,\,i=1,2,\dots,n$ with $\partial_i(x_j)=\delta_{ij}$. Let $\partial$
 be a free 
derivative on $G(K)$. Then $\partial=\sum_{i=1}^nA_i\partial_i$ as 
a free derivative on $F$, where $A_i\in {\mathbb Z}[t^{\pm1}]$,  and it has to 
satisfy the relation
$$\partial(x_j)=t^{\epsilon}\partial(x_k)+(1-t^{\epsilon})\partial(x_i)$$
for each Wirtinger relation $x_j=x_i^{\epsilon}x_kx_i^{-\epsilon}$. Thus the  
${\mathbb Z}[t^{\pm1}]$-module of free derivatives on $G(K)$ can be 
thought of as generated by  
the symbols $A_i,\,i=1,2\dots,n$ and subject to the relation
$$A_j=t^{\epsilon}A_k+(1-t^{\epsilon})A_i$$
for each Wirtinger relation $x_j=x_i^{\epsilon}x_kx_i^{-\epsilon}$.

We define an $n\times n$ matrix $\tilde\B$ as follows. The $j$-th row of 
$\tilde\B$ has at most  
two non-zero entries: for each relation $A_j=t^{\epsilon}A_k+
(1-t^{\epsilon})A_i$, when $k\neq i$, the $(j,k)$-entry is $t^{\epsilon}$ 
and the $(j,i)$-entry is $1-t^{\epsilon}$; 
when $k=i$, the only non-zero entry is the $(j,k)$-entry, which is equal to 1. 

Let $\B$ be the $(n-1)\times(n-1)$ matrix obtained from $\tilde\B$ by deleting
the first
row and the first column. Then 
$\text{det}(I-\B)$ is the Alexander 
polynomial of the knot $K$ (recall that the Alexander polynomial of 
a knot is only defined up to powers of $t$). 
In fact, this is always true no matter which 
$j$-th row and column are deleted. 

\subsection{A model of random walk on knot diagrams} 
In our model of random walk on the knot diagram $K$, we take $\{A_1,A_2,\dots,A
_n\}$ to be the 
space of states. The transition matrix is simply $\tilde\B$. 
This is obviously a 
stochastic matrix
in the sense that the entries in each row add up to 1. In the case when all 
crossings of $K$ are 
positive ($K$ is a positive knot diagram), we get a genuine Markov chain 
for each $t\in [0,1]$. Otherwise, 
we may have
negative probabilities for negative crossings. 

In this model of random walks on $K$, a path from $A_i$ to $A_j$ is 
a sequence of transitions of states from $A_i$ to $A_j$. 
Each such path is associated with a weight (\lq\lq probability"), which is the 
product of \lq\lq transition
probabilities" along this path. 
Pick a state, say $A_1$, consider paths 
from $A_1$ to itself 
which will not contain $A_1$ at any intermediate stage, i.e. we consider paths 
of first return from $A_1$. 
Equivalently, we may regard $A_1$ as being 
broken into two states $A'_1$ 
and $A''_1$, one initial and one terminal. 
This can be done by breaking 
the arc $x_1$ into two arcs $x'_1$ and $x''_1$ and changing 
the knot $K$ into a 
1-string link $T$.
Then we consider all paths on $T$  
from $A'_1$ (the bottom of $T$) to $A''_1$ (the top of $T$). 

\begin{pro} The summation of weights over 
all paths on $T$ from $A'_1$ to $A''_1$ is equal to 1.
\end{pro} 

\begin{proof}  To calculate 
the sum of weights of all paths from $A'_1$ to $A''_1$ amounts to solve the 
system of
linear equations
$$A_j=t^{\epsilon}A_k+(1-t^{\epsilon})A_i$$
for $A''_1$ with $A'_1=1$ given. We have the unique solution $A''_1=1$. For 
more details of the proof,
see \cite{LTW}.
\end{proof} 

We have the following theorem. 
\begin{thm}
1. Let $K$ be a positive knot diagram with $n$ arcs. 
Then for every pair $(i,j)$, there is an integer $m\leq n$, such that the
$(i,j)$-entry of the matrix $\tilde\B^m$ is positive.
Hence, the Markov chain is irreducible.

2. Let $p^{(k)}_{i,j}$ be the $(i,j)$-th entry of 
$\tilde\B^k$.  For each $t\in [0,1]$ and $i,j$, 
$\sum_{k=1}^{\infty} p^{(k)}_{i,j}=\infty$.
Hence each state is persistent.
\end{thm}

\begin{proof} 1.  This is true because we can travel along the knot 
from any state $A_i$ to $A_j$ in $\leq n$ steps.

2. If $i=j$, by Proposition 2.1, 
if we sum the weights of all the $k$-th return paths for $1\leq k\leq n$, the 
sum is $n$. For $i\neq j$, the sum $\sum_{k=1}^{n} p^{(k)}_{i,j}>n$.
\end{proof} 

Imagine that a ball travels on the knot diagram in the 
direction specified by the orientation of the knot.  
It will make a choice when it comes to an $\epsilon$-crossing from 
the under-crossed segment: it may either jump up with probability 
$1-t^{\epsilon}$ and keep traveling on the over-crossed segment or 
keep traveling with probability $t^{\epsilon}$ on the under-crossed segment.  
This is an intuitive picture of our model of random walk on knot diagrams. 
We will call this model the ``jump-up'' model. There is also a ``dual'' model
of jump-down random walk on knot diagrams. In this model, one needs to
make a choice at the over-crossed segment of a crossing: jump-down or keep 
traveling. There are some delicate connections and differences between these
two models which we will not discuss here. We only notice that the two 
random walk models correspond to different 
choices of base points in the Wirtinger presentation.

\subsection{State sum for the Jones polynomial} 
State sum models on knot diagrams is one of the main tools attained in the 
development of topological quantum field theories.  The 
state model we will use for the Jones polynomial 
is given by Turaev in \cite{T} based on earlier constructions of Jones.  
For this model, we need an $R$-matrix.
The $R$-matrix of $\mathfrak {sl}(2)$ with 
respect to the fundamental 
representation is given as follows (with $\bar{q}=q^{-1}$ and $\bar{R}=
R^{-1}$):
$$\begin{aligned}
&R_{0,0}^{0,0}=R_{1,1}^{1,1}=-q,\,R_{0,1}^{1,0}=R_{1,0}^{0,1}=1,
\,R_{0,1}^{0,1}=\bar{q}-q,\\
&\bar{R}_{0,0}^{0,0}=\bar{R}_{1,1}^{1,1}=-\bar{q},
\,\bar{R}_{0,1}^{1,0}=\bar{R}_{1,0}^{0,1}=1,\,\bar{R}_{1,0}^{1,0}=q-\bar{q},
\end{aligned}$$
and all other entries of the $R$-matrix are zero. 

In this model, we consider the 
1-string link $T$ as a planar graph by 
looking at its projection. A state $s$ is an assignment of 0 or 1 to each 
edge of the graph. For each
vertex (crossing) $v$, if $a,b,c,d$ are edges incident to $v$, define
$$\pi_v(s)=(R^{\epsilon})_{s(a),s(b)}^{s(c),s(d)},$$
where $\epsilon$ is the sign of the crossing $v$, $a,b$ are incoming edges and
$c,d$ are outgoing edges.
A state $s$ is admissible
if $\pi_v(s)\neq 0$ for all vertices $v$, and the initial and terminal edges 
having the same assignments. 
The set of all admissible states will be denoted by  
$\text{adm}(T)$. We have 
$$\text{adm}(T)=\text{adm}_0(T)\amalg\text{adm}_1(T)$$
where $\text{adm}_i(T)$, $i=0,1$, is the set of admissible states with $s=i$ on
the initial and terminal edges of $T$. For each admissible state $s$, define
$$\Pi(s)=\prod_{v:\,\text{vertices}}\pi_v(s).$$

Given a 1-string link diagram $T$, and let 
$K$ be a closure of $T$ to a knot diagram without 
introducing any additional crossings, and a state $s\in\text{adm}_i(T), 
i=0,1$. The state $s$ on $T$ can naturally be extended as a state on the 
knot diagram $K$. There are 
quite a few quantities associated with $T$ or the pair $(T,s)$. 
We will define them here, and these notations will be in force throughout
this paper. Also, we will use dashed lines for edges having the assignment 
0 in the state $s$ and solid lines for edges having assignment 1.  

First we define a modification of diagrams according to a state.
{\em A smoothing of $(T,s)$ or $(K,s)$} is the modification 
of the diagram by smoothing 
the crossings marked as
$${\epsfysize=0.12 truein \epsfbox{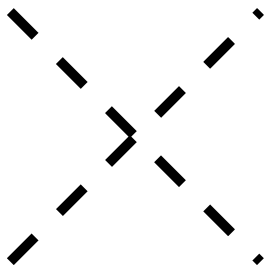}},\,\,
{\epsfysize=0.12 truein \epsfbox{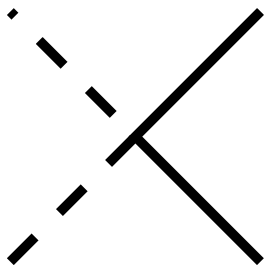}},\,\,
{\epsfysize=0.12 truein \epsfbox{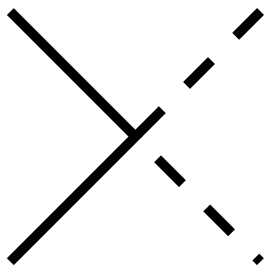}},\,\,
{\epsfysize=0.12 truein \epsfbox{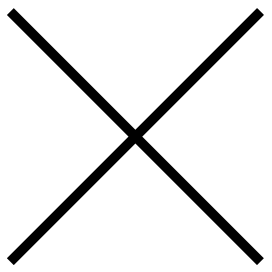}},$$
we get a collection of circles and an arc in the case 
of $(T,s)$, and only circles in the case of $(K,s)$.  
Each circle or arc is marked by $0$ or $1$.

\begin{enumerate}
\item {\em The writhe of $T$:} Denote 
 by $\omega(T)$ the writhe, i.e. the summation of 
signs over all crossings of $T$. 
\item  {$\beta_i(s), i=0,1$:} Denote by $\beta_i(s)$ be the sum of signs 
of crossings whose incident edges all marked by $i$ in $s$.
\item  {\em Rotation numbers}, $\text{rot}(T), \text{rot}_i(K,s), 
\text{rot}_i(T,s)$: Smoothing all 
crossings of $T$, we get a collection of oriented circles in the plane 
(together with an oriented arc), and
$\text{rot}(T)$ is defined to be the sum of rotation numbers 
(Whitney's indices) of these 
circles; For the smoothing of $(T,s)$,
the circles are divided into two collections marked by 
$0$ or $1$ respectively, and $\text{rot}_i(T,s)$ is defined to be the 
sum of rotation numbers of the circles marked by $i$;
The definition of $\text{rot}_i(K,s)$ is similar to that of
$\text{rot}_i(T,s)$, only that the smoothing of $(K,s)$ has one more circle 
then $(T,s)$.
\end{enumerate}

For the Jones polynomial $J(K)$, Turaev's state model gives the following 
formula: 
$$J(K)=(-q^2)^{-\omega(T)}\sum_{s\in\text{adm}(T)}q^{\text{rot}_0(K,s)-
\text{rot}_1(K,s)}\,\Pi(s) .$$
This formula for the Jones polynomial has the value $q+\bar{q}$ on the unknot, 
and the standard variable of the 
Jones polynomial is $t={\bar q}^2$. It is determined by the following 
crossing change formula:
$$\bar t\,J(K_+)- t\,J(K_-)=({\bar t}^{\frac12}-t^{\frac12})\,J(K_0).$$

\noindent{\bf Remark:} This formula is derived from Theorem 5.4 in 
\cite{T}.  The only nontrivial fact is our computation of 
$\int_{D}f$ in the formula which is 
$q^{\text{rot}_0(K,s)-\text{rot}_1(K,s)}$ in our 
notations.  To be more specific, our colors 0,1 correspond to 
the colors 1,2 in \cite{T}, respectively.  Also our conventions for rotation 
numbers are different.  Our convention is that the clockwise oriented 
circle has $\text{rot}=-1$, while the counterclockwise one has $\text{rot}=1$.

Now let us interpret the state sum from the point view of random walks 
on knot diagrams.
First we take a look at the following table:

\medskip
\centerline{\begin{tabular}{r|ccccccc}\hline
$\mathfrak{sl}(2)\,\,\,$ & $\,\,\,
{\epsfysize=0.12 truein \epsfbox{0000.eps}}\,\,\,$ & $\,\,\,
{\epsfysize=0.12 truein \epsfbox{0101.eps}}\,\,\,$ & $\,\,\,
{\epsfysize=0.12 truein \epsfbox{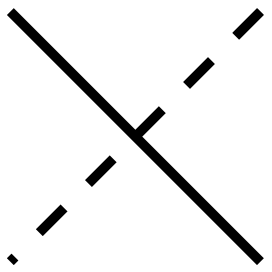}}\,\,\,$ & $\,\,\,
{\epsfysize=0.12 truein \epsfbox{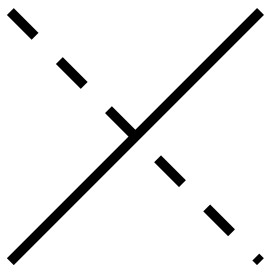}}\,\,\,$ & $\,\,\,
{\epsfysize=0.12 truein \epsfbox{1010.eps}}\,\,\,$ & $\,\,\,
{\epsfysize=0.12 truein \epsfbox{1111.eps}}\,\,\,$
&\,\,\,model\,\,\ \\ \hline
${\epsfysize=0.12 truein \epsfbox{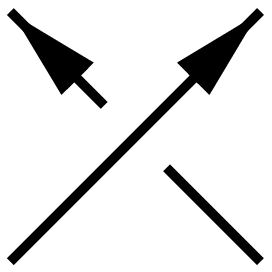}}
\,\,\,$ & $-q$ & $\bar{q}-q$ & 1 & 1 & 0 & $-q$ &\\ \hline
$-\bar q\,{\epsfysize=0.12 truein \epsfbox{plus.eps}}
\,\,\,$ & $1$ & $1-{\bar q}^2$ & ${\bar q}^2\cdot(-q)$ & $(-\bar q)$ & 0 & 
${\bar q}^2q^2$ 
&up\\ \hline
$-\bar q\,{\epsfysize=0.12 truein \epsfbox{plus.eps}}
\,\,\,$ & ${\bar q}^2q^2$ & $1-{\bar q}^2$ & $(-\bar q)$ & ${\bar q}^2\cdot(-q)$ & 0 & $1$ 
&down\\ \hline
${\epsfysize=0.12 truein \epsfbox{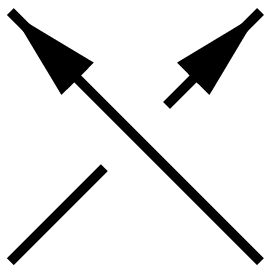}}
\,\,\,$ & $-\bar{q}$ & 0 & 1 & 1 & $q-\bar{q}$ & $-\bar{q}$ &\\ \hline
$-q\,{\epsfysize=0.12 truein \epsfbox{minus.eps}}
\,\,\,$ & $1$ & 0 & $(-q)$ & $q^2\cdot(-\bar q)$ & 
$1-q^2$ & $q^2{\bar q}^2$
&up\\ \hline
$-q\,{\epsfysize=0.12 truein \epsfbox{minus.eps}}
\,\,\,$ & $q^2{\bar q}^2$ & 0 & $q^2\cdot(-\bar q)$ & $(-q)$ & 
$1-q^2$ & $1$
&down\\ \hline
\end{tabular}}
\medskip

Here, as before, a dashed edge has the assignment 0 and a solid edge has 
the assignment 1.
The entry at the row ${\epsfysize=0.12 truein \epsfbox{plus.eps}}$ 
(or $x\,{\epsfysize=0.12 truein \epsfbox{plus.eps}}$) and the column 
${\epsfysize=0.12 truein \epsfbox{0101.eps}}$ is $R_{0,1}^{0,1}$ 
(or $xR_{0,1}^{0,1}$), etc. The last 
column indicates two random walk models  
for this state sum. The two rows marked by ``up'' in the last column 
compare entries of the $xR$ with the weights of 
the jump-up model, and the two rows marked by ``down'' compare entries of
$xR$ with weights of the jump-down model.

Given a state $s\in\text{adm}_0(T)$, think of the edges with assignments $1$ as 
a collection of cycles that several balls traveled in the jump-up model. 
Note that their paths may cross transversely but will not pass through the
same edge twice. Conversely, if we simultaneously have a few balls traveling 
on $T$ avoiding the two open arcs, they do not travel over the same edge but 
may cross transversely, we get a state $s\in\text{adm}_0(T)$ 
by assigning 1 to all the traveled edges, and $0$ otherwise. 
With such a one-one correspondence, for a state $s\in\text{adm}_0(T)$, 
we denote by $W_1^{\circ}(s)$ the 
product of weights of the collection of cycles 
formed by edges marked by 1 as cycles in 
the jump-up model with
$t={\bar q}^2$. 

The case of jump-down model is similar, and it corresponds to 
states in $\text{adm}_1(T)$.  Given such a state $s$,  
the collection of cycles formed by edges marked by 0 are 
thought of as cycles in the jump-down model of random walks and $W_0^{\circ}
(s)$ denotes the product of weights. 

\begin{lm}
In the $\mathfrak{sl}(2)$ state model, we have
$$
\begin{aligned}
&\Pi(s)=(-q)^{\omega(T)}q^{2\beta_1(s)}W_1^{\circ}(s)\qquad 
\text{for $s\in\text{\rm adm}_0(T)$},\\
&\Pi(s)=(-q)^{\omega(T)}q^{2\beta_0(s)}W_0^{\circ}(s)\qquad 
\text{for $s\in\text{\rm adm}_1(T)$}.
\end{aligned}
$$
\end{lm}

\begin{proof} We will show the case $i=0$.  The other case 
is completely similar. 
The factor $(-q)^{\omega(T)}$ comes in since we multiply 
each $R$-matrix entry at an $\epsilon$-crossing by $(-q^{-\epsilon})$.
The term $q^{2\beta_1(s)}$ comes in since we get an 
extra factor $q^{2\epsilon}$ at a solid $\epsilon$-crossing in the 
jump-up model. 
Now using the rows marked by ``up''
in the table above, we need to show the extra multiplicative 
factors of $-q^{\pm1}$ inside $()$ in the columns 
${\epsfysize=0.12 truein \epsfbox{0110.eps}}$ and 
${\epsfysize=0.12 truein \epsfbox{1001.eps}}$ will cancel 
out in the product $\Pi(s)$. 
Notice that after the modification 
of $T$ as we did before, the edges marked by 0 is decomposed into a 
collection of cycles and an arc, having 
transverse intersections with the cycles formed by edges marked 
by 1. The intersections between a
cycle marked by 1 and a cycle or the path marked by 0 can be paired 
up. Consider 
two cases according to whether
such a pair makes a contribution to the linking number. In both cases, we see 
that the extra multiplicative
factors of $- q^{\pm1}$ cancel out.
\end{proof}

Denote
$$
\begin{aligned}
\,&\int_0^0(T)=(-q^2)^{-\omega(T)}\sum_{s\in\text{adm}_0(T)}
q^{\text{rot}_0(T,s)-\text{rot}_1(T,s)}\,\Pi(s),\\
&\int^1_1(T)=(-q^2)^{-\omega(T)}\sum_{s\in\text{adm}_1(T)}
q^{\text{rot}_0(T,s)-\text{rot}_1(T,s)}\,\Pi(s).
\end{aligned}
$$

\begin{lm} We have $\int^0_0(T)=\int^1_1(T)$ and $J(K)=(q+\bar q)\int^0_0(T)=
(q+\bar q)\int^1_1(T).$
\end{lm}

\begin{proof} There are two ways to close up $T$, both giving the same knot $K$.
Thus, we have
$$J(K)=q\int^0_0(T)+\bar q \int^1_1(T)=\bar q\int^0_0(T)+q\int^1_1(T)$$
and this implies the conclusions of the lemma. 
\end{proof}

\subsection{Toward a relationship between Jones polynomial and zeta functions}
Various kinds of zeta functions are basically all about counting 
of cycles. We may also express the Jones polynomial in terms of counting 
``simple families of cycles'' with weights in our model of random walk 
on a $1$-string link $T$. 

Combining previous lemmas, we get the following formula for the 
Jones polynomial.

\begin{lm}  Let $K$ be the closure of a 1-string link $T$,
$$\begin{aligned}
J(K)
&=(q+\bar q)\,q^{-\o(T)+\text{\rm rot}(T)}\sum_{s\in\text{\rm adm}_0(T)}
q^{2(\beta_1(s)-\text{\em rot}_1(T,s))}W_1^{\circ}(s)\\
&=(q+\bar q)\,q^{-\o(T)-\text{\rm rot}(T)}\sum_{s\in\text{\rm adm}_1(T)}
q^{2(\beta_0(s)+\text{\rm rot}_0(T,s))}W_0^{\circ}(s).
\end{aligned}
$$
\end{lm}

\begin{proof} It is not hard to see that $\text{\rm rot}_0(T,s)+
\text{\rm rot}_1(T,s)$ is independent of the state $s$. It is equal to the 
sum of rotation 
numbers of circles obtained by smoothing all crossings of $T$, i.e. the 
rotation number $\text{rot}(T)$ of $T$ by definition.
\end{proof}

To see how the Jones polynomial is related to the Alexander polynomial,
let us describe an expansion of the inverse of the Alexander polynomial.
Consider all cycles in  
our model of random walk which avoid the first arc $A_1$ on the knot diagram. 
Let $\mathcal{Q}$ be the set
of all such cycles which are primitive, i.e. they are 
not powers of any other cycles. 
Recall that $\text{det}\,(I-\mathcal{B})$ is, up to 
a factor of a power of $t$,
the Alexander polynomial of the
knot in question. 
Given a cycle $c$, we will use $W(c)$ to denote its weight.  Then 
$$(\text{det}\,(I-\mathcal{B}))^{-1}=
\prod_{c\in\mathcal{Q}}(1-W(c))^{-1}=1+\sum_{k=1}^\infty\,\sum_{(c_1,\dots,
c_k)\in\mathcal{Q}^k}\,
W(c_1)\cdots W(c_k).$$
This is the Foata-Zeilberger formula we mentioned in the introduction. For
the convenience of readers, an exposition of this formula will be given in 
Section 4.

A $k$-tuple of cycles in ${\mathcal{Q}}^k$ is called {\it simple} if 
no edges are shared by cycles in this $k$-tuple. Let $\mathcal{Q}^t$ 
be the set of all simple $k$-tuples of cycles, for $k=1,2,\dots$.
Given $c\in {\mathcal{Q}}^t$, let $\beta_1(c)$ be the number of crossings 
in $c$, and $\text{rot}(c)$ be the rotation number of $c$.  Note they are 
the same as the $\beta_1$ and $\text{rot}_1$ of the corresponding 
state in $\text{adm}_0(T)$. Finally, in order to have a one-one correspondence
between $\text{adm}_0(T)$ and $\mathcal{Q}^t$, we have to modify $T$ slightly.
The simplest way is to add a negative kink with rotation number $-1$ to the
bottom of $T$ and a positive kink with rotation number 1 to the top of $T$.

We will denote by ${\mathcal Q}_*^t$ the set of all simple $k$-tuples of
cycles in the jump-down model. 

\begin{thm}\label{jexp} With the 1-string link $T$ appropriately chosen as
described above, we have 
$$\begin{aligned}
\frac{J(K)}{q+\bar q}
&=t^{\frac{\o(T)-\text{\rm rot}(T)}{2}}(1+\sum_{c\in{\mathcal{Q}}^t}
t^{\,\text{\em rot}(c)-\beta_1(c)}W(c))\\
&=t^{\frac{\o(T)+\text{\rm rot}(T)}{2}}(1+\sum_{c\in{\mathcal Q}_*^t}
{\bar t}^{(\,\text{\em rot}(c)+\beta_0(c))}W(c)).
\end{aligned}
$$
\end{thm}

Comparing with the expansion of the Alexander polynomial, we see that the 
Jones polynomial $J(K)$ uses the
summands $W(c_1)\cdots W(c_k)$ where no edges are repeated in the collection 
of cycles $c_1,\dots,c_k$. 
A simple idea is that
collections of cycles with repeated edges 
in the expansion of the Alexander polynomial might be lifted to 
collection of simple cycles on the cabling of $T$. 
This idea is realized in Theorem~\ref{mm}.
In the next section, we will first generalize our discussion 
about the Jones polynomial to the colored Jones polynomial. 
 
\section{Limit of the colored Jones polynomial}

\subsection{State sum for the colored Jones polynomial}
The set of finite dimensional irreducible representations of 
$\mathfrak{sl}(2)$ (or rather, the quantum group $U_q\mathfrak{sl}(2)$)
can be listed as $V_1,V_2,V_3,\dots ,$
where $V_d$ is $d$-dimensional. The fundamental representation is $V_2$, which 
is the one used to construct the 
Jones polynomial $J(K)$. Other representations can also be used to construct 
knot polynomials. The knot polynomial
obtained by \lq\lq coloring" the (zero framed) knot $K$ with the 
irreducible representation $V_d$ 
is called the colored Jones polynomial $J(K,V^d)$ 
\cite{RT}. We have $J(K,V_1)=1$ and $J(K,V_2)=J(K)$. And if $K$ is the unknot,
$$J(K,V_d)=[d]=\frac{q^d-\bar{q}^d}{q-\bar{q}}.$$

We may also color $K$ by non-irreducible representations, for example, by
$V_2^{\otimes d}$. Such a colored Jones polynomial can be interpreted in 
two ways:
\begin{enumerate}
\item Assume that $K$ has zero framing, let $K^d$ be the link obtained by 
replacing $K$ with $d$ parallel copies
(this is the zero framing cabling operation), then $J(K,V_2^{\otimes d})=
J(K^d).$
\item We have the following relation in the representation ring of $\mathfrak
{sl}(2)$: $V_2\otimes V_d=V_{d+1}\oplus
V_{d-1}$. Thus, $V_2^{\otimes d}$ is a linear combination 
of the irreducible modules $V_{d+1}$, 
$V_{d-1}$, $V_{d-3}$, ... and 
$J(K,V_2^{\otimes d})$ is the same linear combination of $J(K,V_{d+1})$, 
$J(K,V_{d-1})$, $J(K,V_{d-3})$, ....
\end{enumerate}

These two interpretations can be used to establish a precise relation between 
the colored Jones polynomials and
the cablings of the Jones polynomial. We quote from \cite{KM} such a relation 
in the case considered here:
$$J(K,V_{d+1})=\sum_{j=0}^{d/2}\,(-1)^j\,\begin{pmatrix} d-j\\ j\end{pmatrix}
\, J(K^{d-2j})\,.$$ 

The decomposition $V_2\otimes V_d=V_{d+1}\oplus V_{d-1}$ can be given 
explicitly in terms of the standard bases of these 
irreducible representations \cite{KR}. 
Suppose the 
standard basis of $V_2$ is $\{e_0,e_1\}$, and the standard basis of $V_{d+1}$
is $\{f_0,f_1,\dots,f_d\}$, then we have 
$$\begin{aligned}
\,&f_0=a\cdot e_0\otimes e_0\otimes\cdots\otimes e_0\in V_2^{\otimes d},\\
&f_d=b\cdot e_1\otimes e_1\otimes\cdots\otimes e_1\in V_2^{\otimes d},
\end{aligned}$$
where $a,b$ are products of $q$-analogue of Clebsch-Gordan coefficients 
\cite{KR}.

For a 1-string link $T$, if it is colored by $V_{d+1}$, we get an invariant
$F(T)$ which is a $U_q\mathfrak{sl}(2)$-morphism of $V_{d+1}$. Since $V_{d+1}$
is an irreducible $U_q\mathfrak{sl}(2)$-module, we have 
$$F(T)(f_i)=\lambda\,f_i,\qquad i=0,1,\dots,d.$$
Furthermore, let $K$ be the closure of $T$, then
$$J(K,V_{d+1})=[d+1]\cdot\lambda.$$

On the other hand, if we color $T$ by $V_2^{\otimes d}$, we may write the 
induced
$U_q\mathfrak{sl}(2)$-morphism $F(T)$ on $V_2^{\otimes d}$ as follows:
$$F(T)(e_{i_1}\otimes\cdots\otimes e_{i_d})=
\sum_{j_1,\dots,j_d}\,\int_{i_1\cdots i_d}
^{j_1\cdots j_d}(T)\,e_{j_1}\otimes\cdots\otimes e_{j_d}.$$
Thus, the following lemma holds, which generalizes Lemma 2.4. 

\begin{lm} We have $\int_{0\cdots0}^{0\cdots0}(T)=
\int_{1\cdots1}^{1\cdots1}(T)=\lambda$ 
and 
$J(K,V_{d+1})=[d+1]
\,\int_{0\cdots0}^{0\cdots0}(T)=[d+1]\,\int_{1\cdots1}^{1\cdots1}(T)$.
\end{lm}

We now can extend Theorem~\ref{jexp} 
to $J(K,V_{d+1})$. Notice that we assume the
writhe $\o(T)=0$ and $T^d$ is the zero-framing $d$-cabling of $T$. We denote
by $\text{adm}_0(T^d)$ the set of admissible states on $T^d$ which assign 0 to 
all the top and bottom edges. The notation $\text{adm}_1(T^d)$ has the 
obvious  
meaning.

\begin{lm}\label{color} With the notations as above, we have 
$$\begin{aligned}
J(K,V_{d+1})
&=[d+1]\,q^{\text{\rm rot}(T^d)}\sum_{s\in\text{\rm adm}_0(T^d)}
q^{2(\beta_1(s)-\text{\em rot}_1(T^d,s))}W_1^{\circ}(s)\\
&=[d+1]\,{\bar q}^{\text{\rm rot}(T^d)}\sum_{s\in\text{\rm adm}_1(T^d)}
q^{2(\beta_0(s)+\text{\rm rot}_0(T^d,s))}W_0^{\circ}(s).
\end{aligned}
$$
\end{lm}

\begin{proof} Applying Turaev's state model to the tangle $T^d$, we get
$$\int_{i_1\cdots i_d}^{j_1\cdots j_d}(T)=(-q^2)^{w(T^d)}\sum_{s\in\text{adm}_*(T^d)}
\,q^{\text{rot}_0(T^d,s)-\text{rot}_1(T^d,s)}\Pi(s)$$
where $\text{adm}_*(T^d)$ is the set of admissible states on $T^d$ such that
the bottom edges are assigned with $i_1,\dots,i_d$ and top edges with
$j_1,\dots,j_d$, respectively. Then we can translate this expression for  
$\int_{i_1\cdots i_d}^{j_1\cdots j_d}(T)$ into the form appeared in 
Proposition~\ref{color} as we did in Section 2.4.
\end{proof}

Now the corresponding formula 
for $J(K,V_{d+1})$ of  Theorem~\ref{jexp}  is 
obtained from Theorem~\ref{jexp} by replacing 
$q+\bar q=[2]$ by $[d+1]$. 

\subsection{Computation of the limit}  
In this section, we prove our main theorem which calculates 
the limit of the renormalized colored Jones polynomials 
when the color
parameter tends to infinity and the weight parameter tends to 1.

\begin{thm}\label{mm} 
Let $T$ be a 0-framed 1-string link, modified appropriately as in the 
Theorem 2.6, and $K$ be the closure of $T$. Denote by $\mathcal Q$ 
(${\mathcal Q}_*$)
the set of primitive cycles in the jump-up (jump-down) model of random walk 
on $T$ with $t=e^{-2h}$. Then
$$\begin{aligned}
\lim_{d\rightarrow\infty}\frac{J(K,V_{d+1})(e^{\frac{h}{d}})}{[d+1]} 
&={\bar t}^{\,\frac{\text{\rm rot}(T)}2}\,\left(1+\sum_{(c_1,c_2,\dots,c_k)\in
{\mathcal Q}^k}\,W(c_1)W(c_2)\cdots W(c_k)\right)\\
&={t}^{\,\frac{\text{\rm rot}(T)}2}\,\left(1+\sum_{(c_1,c_2,\dots,c_k)\in
{\mathcal Q}_*^k}\,W(c_1)W(c_2)\cdots W(c_k)\right).
\end{aligned}
$$
\end{thm}

\begin{proof}
Using the expansion of the colored Jones polynomials, it 
suffices to show that the weight of cycles 
$(c_1, c_2, \cdots , c_k)$ on $T$ in the right-handed side 
is the limit of some cycles on $T^d$ for large $d$. 
Let us compare the two jump-up models of random walks on $T$ 
and $T^d$ with $t=e^{-2h}$, 
and with $t=e^{-\frac{2h}{d}}$, respectively. 

Consider first a simple cycle $c$ 
on $T$. Recall that this is a cycle on $T$ with no edges repeated. 
There are many ways to lift $c$ to become a simple cycle $\tilde c$
on $T^d$. The reason for this multiplicity is that for each jump-up on $c$,
we can choose one of the $d$ over-crossed segments to jump up on $T^d$. In fact,
if there are $m$ jump-ups on $c$, there will be $m^d$ lifts $\tilde c$ on 
$T^d$. We need to calculate $\sum\,W(\tilde c)$, a sum over all liftings of 
$c$. For a jump-up at a positive crossing on $c$, we get a (multiplicative) 
contribution 
$1-e^{-2h}$ to $W(c)$. The corresponding contribution to $\sum\,W(\tilde c)$
is a multiplicative factor
$$(1-e^{-\frac{2h}{d}})(1+e^{-\frac{2h}{d}}+e^{-\frac{4h}{d}}+\cdots+
e^{-\frac{2(d-1)h}{d}})=1-e^{-2h}.$$
Also, passing through an under-crossing on $c$ contributes $e^{-2h}$ to $W(c)$
and the corresponding contribution of $\tilde c$ is 
$$(e^{-\frac{2h}{d}})^d=e^{-2h}.$$
Thus we have 
$$\sum\,W(\tilde c)=W(c).$$

Obviously,
$\beta_1(\tilde c)$ and $\text{rot}_1(T^d,\tilde c)$ depend only on $c$. 
We also notice that
$\text{rot}(T^d)=d\,\text{rot}(T)$. Thus, 
$$\begin{aligned}
\lim_{d\rightarrow\infty}&\,(e^{\frac{2h}{d}})^{\,\frac{\text{rot}(T^d)}{2}}
\sum_{\tilde c}
(e^{-\frac{2h}{d}})^{\text{rot}_1(T^d,\tilde c)-\beta_1(\tilde c)}
W(\tilde c)\\
&=(e^{2h})^{\,\frac{\text{\rm rot}(T)}2}W(c).
\end{aligned}
$$ 

Notice that the same argument holds true for a simple collection of cycles 
on $T$.

In general, given a non-simple collection of cycles $c$ on $T$, we decorate 
each edge by an integer which is the number of times $c$ traveling over  
that edge. There are only finitely many collections of cycles on $T$ with a 
fixed decoration. For $d$ sufficiently large, we can lift $c$ to a simple 
collection of cycles on $T^d$. To get such a lifting, we will not have the
freedom of jumping up onto any of the $d$ over-crossed segments at a crossing.
A particular jump-up at a crossing $X$ on $T$ has at most $d$ liftings.
 For some other 
jump-up onto the segment going over $X$, we have to avoid 
the over-crossed segments jumped onto previously. There are at most $d$ possible collisions 
for the liftings of these two jump-ups. Since
$$\lim_{d\rightarrow\infty}\,(1-e^{\pm\frac{2h}{d}})(1-e^{\pm\frac{2h}{d}})d=0,$$ 
we conclude that in the limit when $d\rightarrow\infty$, 
the sum of weight of all
non-simple liftings of $c$ is zero. We may just do our calculation
as if there are only simple liftings. Thus,
the same calculation as we did before leads to
$$\lim_{d\rightarrow\infty}\sum_{\text{$\tilde c$ simple}}\,W(\tilde c)
=W(c).$$
 
Finally, 
$\beta_1(\tilde c)$ and 
$\text{rot}_1(T^d,\tilde c)$ are bounded by quantities depending only on $c$. 
Thus, we get
$$\lim_{d\rightarrow\infty}\frac{J(K,V_{d+1})(e^{\frac{h}{d}})}{[d+1]} 
={\bar t}^{\frac{\text{\rm rot}(T)}2}\,\left(1+\sum_{(c_1,c_2,\dots,c_k)\in
{\mathcal Q}^k}\,W(c_1)W(c_2)\cdots W(c_k)\right).
$$

This finishes the proof of Theorem~\ref{mm}.
\end{proof}

\section{Ihara-Selberg zeta function and Melvin-Morton conjecture}

\subsection{Lyndon words and the Foata-Zeilberger formula} 
Let us recall the notion of Lyndon words and some results in \cite{FZ}. For 
references to quoted results in this section, see \cite{FZ}.

Given a finite nonempty set $X$ whose elements are totally ordered, 
we consider the monoid $X^*$ generated by $X$. Let $<$ be the lexicographic 
order on $X^*$ derived from the 
total order on $X$. A {\it Lyndon word} is defined to be a nonempty word in 
$X^*$ which is prime, i.e. not the 
power of any other word, and is minimal in the class of its cyclic 
rearrangements. Let $L$ denote the set of all
Lyndon words. The following result is due to Lyndon.

\begin{lm} Each nonempty word $w\in X^*$ can be uniquely written as a 
non-increasing juxtaposition of Lyndon words:
$$w=l_1l_2\cdots l_m,\qquad l_i\in L,\,\,\, l_1\geq l_2\geq\cdots\geq l_m.$$
\end{lm}

Let $X$ be a finite set. Let $\mathcal B$ be a square matrix whose entries 
$b(x,x')$ ($x,x'\in X$) form a set of commuting 
variables. For each Lyndon word $l\in L$, we associate with it a variable 
denoted by $[l]$. These variables
$[l]$, $l\in L$, are assumed to be all distinct and commute with each other.

Given a word $w=x_1x_2\cdots x_k$ in $X^*$, define
$$\beta_{\text{circ}}(w)=b(x_1,x_2)b(x_2,x_3)\cdots b(x_{k-1},x_k)b(x_k,x_1)$$
and $\beta(w)=1$ if $w$ is the empty word. Notice that all the words in the 
same cyclic rearrangement class have the 
same $\beta_{\text{circ}}$-image. Also define
$$\beta([l])=\beta_{\text{circ}}(l)$$
for $l\in L$.

Now form the $\mathbb Z$-algebras of formal power series in the variables 
$[l]$ and $b(x,x')$ respectively. 
Extend $\beta$ to a continuous homomorphism between these two 
$\mathbb Z$-algebras. It makes sense
to consider the product
$$\Lambda=\prod_{l\in L}(1-[l])$$
as well as its inverse $\Lambda^{-1}$. We have
$$\beta(\Lambda)=\prod_{l\in L}(1-\beta[l])$$
and 
$$\beta(\Lambda^{-1})=(\beta(\Lambda))^{-1}.$$

For a nonempty word $w\in X^*$, let it be written as in Lemma 4.1. Then define
$$\beta_{\text{dec}}(w)=\beta_{\text{circ}}(l_1)\beta_{\text{circ}}(l_2)
\cdots\beta_{\text{circ}}(l_k).$$
If $w$ is empty, $\beta_{\text{dec}}(w)=1$. Finally, define
$$\beta_{\text{dec}}(X^*)=\sum_{w\in X^*}\,\beta_{\text{dec}}(w).$$

The following theorem of Foata and Zeilberger is what we need.

\begin{thm}\label{FZfor} {\em (Foata-Zeilberger formula)} 
$\beta(\Lambda^{-1})=\beta_{\text{\rm dec}}(X^*)=(\text{\rm det}\,
(I-\mathcal{B}))^{-1}.$
\end{thm}

This is a generalization of the Bowen-Lanford formula \cite{BL}, 
which comes directly
from the identity $\text{det}(e^{A})=e^{\text{tr}A}$ for a matrix $A$.
 
\subsection{The Ihara-Selberg zeta function of a graph} The 
Foata-Zeilberger formula in
Theorem 3.2 is used in \cite{FZ} to derive one of
Bass' evaluations of the Ihara-Selberg zeta function for a graph \cite{B}. For 
the reader's convenience, let 
us first recall Ihara's formulation of the zeta function in the original 
setting of Selberg (see \cite{B}).

Let $\Gamma<PSL_2(\mathbb{R})$ be a uniform lattice (= discrete cocompact 
subgroup). An element $g\in\Gamma$ is hyperbolic if
$$l(g)={\text{min}}\{d\,(gx,x)\,;\,x\in\mathbb{R}^2_+\}>0\qquad \text{($d=$ 
Poincar\'e metric)}.$$

Let $\mathcal{P}$ be the set of $\Gamma$-conjugacy classes of primitive 
hyperbolic elements
in $\Gamma$, then the Ihara-Selberg zeta function is
$$Z(s)=\prod_{g\in\mathcal{P}}(1-u^{l(g)})^{-1},\qquad u=e^{-s}.$$

Let $G$ be an directed graph with the set of edges $E(G)=\{e_1,e_2,\dots,
e_n\}$.
Let $S$ be an $n\times n$ matrix
whose $(i,j)$-entry is equal to 1 if the terminal point of $e_i$ is the same 
as the initial point of $e_j$, and 
0 otherwise. On $G$, we may consider primitive cycles, which are oriented 
cycles formed by directed edges in the
usual sense and which are not powers of some other cycles. Let $\mathcal{C}$ 
be the set of
primitive cycles on $G$. The Ihara-Selberg zeta function of $G$ is 
$$Z_G(u)=\prod_{c\in\mathcal{C}}(1-u^{|c|})^{-1},$$
where $|c|$ is the length of the cycle $c$ (= the number of edges in $c$). 
The Foata-Zeilberger formula implies
$$Z_G(u)=(\text{det}\,(I-uS))^{-1}.$$

If $G$ is an undirected graph, in \cite{B}, Bass transformed $G$ into an 
directed graph $G'$ by 
giving each edge of $G$ two different orientations and thinking of them as 
different directed edges. To study 
primitive, reduced cycles on $G$, where \lq\lq reduced" means that an edge 
will not be traveled twice successively, Bass 
looked at the matrix $T=S-J$, where $S$ is the matrix we defined in the 
previous paragraph for $G'$ and $J$ 
is the matrix whose $(i,j)$-entry is 1 if the $i$-th and $j$-th edges of $G'$ 
come from the same edge of $G$, 
and 0 otherwise. Now let $\mathcal{R}$ be the set of primitive, reduced 
cycles on $G$, define
$$Z_G(u)=\prod_{c\in\mathcal{R}}\,(1-u^{|c|})^{-1}.$$
One of Bass' evaluations of $Z_G(u)$, which is now a consequence of 
the Foata-Zeilberger formula, is
$$Z_G(u)= (\text{det}\,(I-uT))^{-1}.$$
 
The Foata-Zeilberger formula is general enough so that we may apply it to 
Markov processes with a finite set of 
states. A cycle now will be a sequence of transitions of states from  
and back to a given one. In particular, in our model of random walk on a knot 
diagram discussed in Sections 2.1 
and 2.2, we have the set of states $\{A_1,A_2,\dots,A_n\}$ and the transition 
matrix $\tilde{\mathcal{B}}$. 
This case is degenerate since $\text{det}\,(I-\tilde{\mathcal{B}})=0$. 
Nevertheless, we may consider all cycles in  
our model of random walk which avoid the first arc $A_1$ on the knot diagram. 
Let $\mathcal{Q}$ be the set
of all such cycles which are primitive, then the 
Foata-Zeilberger formula implies
$$\prod_{c\in\mathcal{Q}}\,(1-W(c))^{-1}=(\text{det}\,(I-\mathcal{B}))^{-1},$$
where $W(c)$ is the weight of the cycle $c$ and $\mathcal{B}$ is obtained from 
$\tilde{\mathcal{B}}$ by deleting 
the first row and column. Notice that $\text{det}\,(I-\mathcal{B})$ is, up to 
a factor of a power of $t$,
the Alexander polynomial of the
knot in question. So we see that the inverse of the Alexander polynomial is an 
Ihara-Selberg type zeta function.

We have
$$\prod_{c\in\mathcal{Q}}(1-W(c))^{-1}=1+\sum_{k=1}^\infty\,\sum_{(c_1,\dots,
c_k)\in\mathcal{Q}^k}\,
W(c_1)\cdots W(c_k).$$
Hence we obtain the following expansion of the inverse of the 
Alexander polynomial:
\begin{thm}\label{aexp}
$$(\text{det}\,(I-\mathcal{B}))^{-1}=
\prod_{c\in\mathcal{Q}}(1-W(c))^{-1}=1+\sum_{k=1}^\infty\,\sum_{(c_1,\dots,
c_k)\in\mathcal{Q}^k}\,
W(c_1)\cdots W(c_k).$$
\end{thm}

\subsection{Melvin-Morton function and Melvin-Morton Conjecture}

In \cite{MM}, Melvin and Morton 
studied the dependence of the colored Jones polynomial on the
``color'' (that is the dimension $d$). They observed that
$$\frac{J(K,V_{d+1})(e^h)}{[d+1]}=\sum_{m\geq0,\,j\leq m}\,a_{jm}(K)d^jh^m.$$
Furthermore, Melvin and Morton conjectured that the function (which will be 
called the {\it Melvin-Morton function})
$$M(K)(h)=\sum_{m\geq0}\,a_{mm}(K)h^m$$
is the inverse of the Alexander polynomial.

Rozansky then was able to give a proof of this conjecture, on the level of
rigor of physics, based essentially on calculating the limit
$$\lim_{d\rightarrow\infty}\frac{J(K,V_{d+1})(e^{\frac{h}{d}})}{[d+1]}$$
and the known relationship between the semi-classical limit of Witten's
Chern-Simons path integral and the Ray-Singer torsion. Rozansky's work
went beyond the particular simple Lie algebra $\mathfrak{sl}(2)$ and extended 
the Melvin-Morton conjecture to its full generality.

The first rigorous proof of the Melvin-Morton conjecture was given by
Bar-Natan and Garoufalidis \cite{BG}. Their proof used the full power of 
the theory of finite type knot invariants, together with some quite complicated 
combinatorial arguments. Later, Vaintrob and others simplified the 
combinatorial arguments of Bar-Natan and Garoufalidis (see, for example, \cite{V}). 

The Melvin-Morton conjecture can be deduced now as follows. 
By Theorem~\ref{mm} and Theorem~\ref{aexp},
$$\lim_{d\rightarrow\infty}\frac{J(K,V_{d+1})(e^{\frac{h}{d}})}{[d+1]} 
=\frac{{\bar t}^{\frac{\text{\rm rot}(T)}2}}{\text{\rm det}(I-{\mathcal B})}.
$$
On the other hand, it is easy to see that 
$$\lim_{d\rightarrow\infty}\frac{J(K,V_{d+1})(e^{\frac{h}{d}})}{[d+1]}
=M(K)(h).$$
Hence the Melvin-Morton conjecture follows:

\begin{thm}
For any knot $K$ which is the closure of a 0-framed 1-string link $T$,
$$M(K)(h)=
\frac{{\bar t}^{\,\frac{\text{\rm rot}(T)}2}}{\text{\rm det}(I-{\mathcal B})},
$$ where $t=e^{-2h}$. 
\end{thm}

Note the right side is the inverse of the symmetric Alexander 
polynomial of $K$ when the 1-string link $T$ is chosen appropriately as in
Theorem 2.6.

\noindent{\bf Remark:} In Theorem 3.3, we are actually calculating the 
limit of the
partition function $\int_{0\cdots0}^{0\cdots0}(T)$ 
with a fixed boundary condition. 
This is rather like the calculation in statistical mechanics (e.g 
the limit of the Ising model). In statistic mechanics,  
the discontinuities of the limiting function are
related with phase transitions. Thus, it might make sense to ask 
whether the zeros of the Alexander polynomial are of any significance and 
could be ``observed''.


\begin{thebibliography}{99999}

\bibitem{B} H. Bass, {\it The Ihara-Selberg zeta function of a tree lattice,} 
Internat. J. Math. 
{\bf 3}(1992). 717--797. 

\bibitem{BG} D. Bar-Natan and S. Garoufalidis, {\it On the 
Melvin-Morton-Rozansky conjecture}, Invent. Math. {\bf 1}(1996), 103-133. 

\bibitem{BL} R. Bowen and O. Lanford, {\it Zeta functions of restrictions 
of the shift transformations, } Proc. Symp. Pure. Math., {\bf 14}(1995),
43-49.

\bibitem{FZ} D. Foata and D. Zeilberger, {\it A combinatorial proof of Bass' 
evaluations of
the Ihara-Selberg zeta function for graphs}, Trans. AMS, to appear. 
E-print: {\tt math.CO/9806037}.

\bibitem{KM} R. Kirby and P. Melvin, {\it The 3-manifold invariants of 
Witten and 
Reshetikhin-Turaev for $\mathfrak{sl}(2)$,} Invent. Math., {\bf 105}(1991), 
473-545. 

\bibitem{KR} A. N. Kirillov and N. Yu. Reshetikhin, {\it Representations
of the algebra $U_q\mathfrak{sl}(2)$, $q$-orthogonal polynomials and 
invariants of links}, Infinite Dimensional Lie Algebras and Groups
(ed. V. G. Kac), World Scientific, 1989, pp. 285--339.

\bibitem{LTW} X.-S. Lin, F. Tian and Z. Wang, {\it Burau representation and 
random walk on string links,}
Pacif. J. Math., {\bf 182}(1998), 289-302. 

\bibitem{MM} P. Melvin and H. Morton, {\it The coloured Jones function}, 
Commun. Math. Phys., {\bf 169}(1995), 501-520.

\bibitem{RT} N. Yu. Reshetikhin and V. G. Turaev, {\it Ribbon graphs and 
their invariants derived from 
quantum groups,} Comm. Math. Phys. {\bf 127}(1990), 1--26.

\bibitem{RO} L. Rozansky, {\it A contribution of the trivial connection 
to the Jones polynomial and Witten's invariant of 3-manifolds,} 
 Comm. Math. Phys. {\bf 175}(1996), 275-318.

\bibitem{T} V. G. Turaev, {\it The Yang-Baxter equation and invariants of 
links,} Invent. Math. 
{\bf 92}(1988), 527-553.

\bibitem{V}A. Vaintrob, {\it Melvin-Morton conjecture and primitive
Feynman diagrams,} Internat. J. Math. {\bf 8} (1997), 537--553.


\end{thebibliography}
\end{document}